\documentclass[a4paper,10pt]{article}
\usepackage[latin1]{inputenc}
\usepackage{amsmath, amsthm, amssymb}
\usepackage{fullpage}

\newtheorem{theo}{Theorem}
\newtheorem{pr}{Proposition}
\newtheorem{cor}{Corollary}
\newtheorem{lem}{Lemma}

\newcommand{\CQFD}{\hfill $\square$}

\newcommand{\cE}{{\cal E}}
\newcommand{\cF}{{\cal F}}
\newcommand{\cG}{{\cal G}}

\newcommand{\cM}{{\cal M}}
\newcommand{\cN}{{\cal N}}

\newcommand{\cP}{{\cal P}}

\newcommand{\cS}{{\cal S}}
\newcommand{\cT}{{\cal T}}

\newfont{\msbm}{msbm10 scaled\magstep1}
\newfont{\msbms}{msbm7 scaled\magstep1} 

\newcommand{\bbE}{\mbox{$\mbox{\msbm E}$}}

\newcommand{\bbN}{\mbox{$\mbox{\msbm N}$}}

\newcommand{\bbP}{\mbox{$\mbox{\msbm P}$}}

\newcommand{\bbR}{\mbox{$\mbox{\msbm R}$}}

\title{Phenotypic diversity and population growth in fluctuating environment: a MBPRE approach.}
\author{ Clément Dombry\thanks{Laboratoire de Math\'ematiques et Applications, T\'el\'eport 2- BP30179, Boulevard Pierre et Marie Curie, 86962 Futuroscope Chasseneuil Cedex, France. clement.dombry@math.univ-poitiers.fr},\ \ Christian Mazza\thanks{D\'epartement de Math\'ematique, Université de Fribourg, Chemin du Mus\'ee 23, CH-1700 Fribourg, Suisse.Christian.Mazza@unifr.ch}\ \ and Vincent Bansaye\thanks{CMAP, Ecole Polytechnique,  Route de Saclay,
91128 Palaiseau Cedex, FRANCE. bansaye@polytechnique.edu} \protect\hspace{1cm}}
\begin{document}

\maketitle

\begin{abstract}
Organisms adapt to fluctuating environments by regulating their dynamics, and by
adjusting their phenotypes to environmental changes. We model population growth
using multitype branching processes in random environments, where the offspring
distribution of some organism having trait $t\in\cT$ in environment
$e\in\cE$ is given by some (fixed) distribution $\Upsilon_{t,e}$ on $\bbN$. Then, the phenotypes
are attributed using a distribution (strategy) $\pi_{t,e}$ on the trait space $\cT$. We look
for the optimal strategy $\pi_{t,e}$, $t\in\cT$, $e\in\cE$ maximizing the net growth rate
or Lyapounov exponent, and characterize the set of optimal strategies. This is considered for
various models of interest in biology: hereditary versus non-hereditary strategies and strategies
involving or not involving a sensing mechanism. Our main results are obtained in the setting of non-hereditary strategies: thanks to a reduction to simple branching processes in random environment, we derive an exact expression for the net growth rate and a characterisation of optimal strategies. We also focus on typical genealogies, that is, we 
consider the problem of finding the typical lineage of a randomly chosen organism.
\end{abstract}
\\

\section{Introduction\label{s.intro}}
Organisms adapt to fluctuating environments by regulating their intrinsic dynamics, and adjust their phenotypes or traits to
the random environment. Observations reveal that most cell populations are heterogeneous, that is, are composed of various phenotypes.
In \cite{HacIwa}, \cite{KuLei} and \cite{Thattai}  basic models are provided to explain this heterogeneity: Time is continuous, and the related stochastic processes can be seen as 
 birth processes  for the different phenotypes  with migration processes between them in fluctuating environments.   We first present the models considered in \cite{KuLei}, and then
turn to our new setting. In \cite{KuLei}, the trait space $\cT$ is finite, and  the environment space $\cE$ is such that $\cE =\cT$, the idea being that trait $i$ has the fastest growth rate in environment $e=i$.  The migration rate matrix $H^{(k)}=(H_{ij}^{(k)})_{i,j\in\cT}$ gives the switching rates $H_{ij}^{(k)}$ from phenotype $j$ to $i$ in environment $e=k$. The changing environment is modeled as an alternating renewal process. The migration rates $H_{ij}^{(k)}$ defines possible strategies to overcome the uncertain future. \cite{KuLei}
considered two basic strategies: {\it stochastic switching} and {\it responsive switching}. In the first case, the organisms decide to switch the phenotypes independently of the running environment, they do not use sensors, so that the migration rates can be written as
$H_{ij}^{(k)}=H_{ij}$. Responsive switching assumes strong use of sensors, at a certain cost, and the extreme way of responding to environmental changes is obtained by assuming migration rates of the form $H_{kj}^{(k)} = H_m$ when $j\ne k$, and $H_{ij}^{(k)} =0$ when
$i\ne k$ and $j\ne i$. In this last situation, organisms opt always for the most favourable phenotype in the running environment, that is for the trait having the largest growth rate. Denoting by $X_i(t)$ the mean number of organisms having phenotype $i$ at time $t$, the average $X(t)$ is solution of the random differential equation $dX(t)/dt = A_{\cE (t)}X(t)$, where $A_{\cE (t)}$ denotes the rate matrix composed of the various growth and migration rates in environment $\cE (t)$ at time $t$. The authors then estimate the related Lyapounov exponent $\ln(N(t))/t$ as $t\to\infty$, where $N(t)=\sum_{i\in \cE}X_i(t)$ is the average total population size. The computations are performed under the assumption that the environment changes slowly, so that, the process has the time to attain an equilibrium during each phase, see \cite{KuLei} for more details. A similar model was also considered by \cite{Thattai} for two phenotypes and two environments. The model was studied using  Monte-Carlo simulations, which indicate that responsive switching is sometimes not the optimal strategy. The average proportion of cells having the fastest growth rate is maximized for strategies using a sensing mechanism allowing migration to the unfavourable state, leading to heterogeneous populations, providing in this way a rationale for population diversification. These models were then studied analytically in \cite{Gander}. Modifications of these models were also considered more recently in \cite{Visco} where catastrophic events are modeled. 

In the present work, we treat similar problems using discrete time multitype branching processes in random environments, 
denoted by MBPRE in what follows,
 which offer alternatives to the above birth with migration processes. The advantage of multitype  branching processes 
 is that one can separate more clearly the birth and migration phases. We can for example treat in this way populations composed of
 organisms where birth occurs in a very precise period. When this is not the case, discrete time modeling still provides
 relevant informations for population growth.
  We will also consider the problem of finding optimal strategies to maximize
 the random Lyapounov exponent.
 Our
 mathematical approach also permits to treat situations where the trait and environment spaces $\cT$ and $\cE$ are continuous: all of our results are valid in very general situations, like the models considered in \cite{HacIwa}. To be useful for scientists not familiar with too advanced mathematics, we first illustrate some
 results and provide examples when both $\cT$ and $\cE$ are finite.
\subsection{Results\label{s.results}}


Assume that $\cT =\{t_1,\cdots,t_q\}$ and that $\cE =\{e_1,\cdots,e_p\}$.
The process of interest  $(Z_n)_{n\ge 0}$ is written as a random vector $Z_n = (Z_n^t)_{t\in\cT}$, where $Z_n^t$
denotes the number of organisms having trait $t\in\cT$ at generation $n$. The transition between generations $n$
and $n+1$ is modeled as a two step procedure. First, for given $t\in\cT$, each of the $i=1,\cdots,Z_n^t$
organisms having trait $t$ gives birth to a random number of descendants given by a random variable
$\xi_{i,n}^{t,\omega_n}$, where the index $\omega_n\in\cE$ models the environmental state at generation $n$. In the second phase, each of the $j=1,\cdots,\xi_{i,n}^{t,\omega_n}$ new individuals
is assigned to a random trait $\tau_{i,j,n}^{t,\omega_n}\in\cT$. As usual in such processes, we assume 
that all of these random variables are independent, with i.i.d. $\xi_{i,n}^{t,\omega_n}$ and
i.i.d. $\tau_{i,j,n}^{t,\omega_n}$. We further assume that the environmental process
$(\omega_n)_{n\ge 0}$ is fixed, and is supposed to be the realization of a stationary 
\footnote{i.e. for all  $i \in \mathbb N$, $(w_0,w_1,...)$ is distributed as $(w_{i},w_{i+1},...)$} 
and ergodic process \footnote{i.e. for all bounded Borel function $f$, $\mathbb P(f(w_0,w_1,...)=f(w_1,w_2,...))=1$ implies  $f(w_0,w_1,...)$ is constant.} taking
values in $\cE$. Let  $\Upsilon_{t,e}$ and $\pi_{t,e}$ be distributions of the $\xi_{i,n}^{t,e}$ and
$\tau_{i,j,n}^{t,e}$ for a trait $t\in\cT$ in environment $e\in\cE$. We will often use the 
first moment $m_{t,e}$ of $\Upsilon_{t,e}$, giving the average number of descendants for an organism
having trait $t$ in environment $e$. Optimal strategies will be considered for fixed distributions 
$\Upsilon_{t,e}$, that is, {\it we will maximize the population growth rate as a function of the
distributions $\pi_{t,e}$}. Sections 3, 4 and 5 consider extinction and optimal growth questions. Section 6
deals with typical genealogies, that is, we consider the question of finding the typical lineage
of an individual chosen randomly at generation $n$. We also illustrate these notions both for finite
spaces $\cT$ and $\cE$ and for continuous ones.

As stated in the previous Section, the authors in \cite{KuLei} distinguish between stochastic and responsive switching. In the first case,
the distributions $\pi_{t,e}$ on $\cT$ depends on $t$ but not on $e$, so that we can write $\pi_{t,e}\equiv \pi_t$. Concerning the second
family of strategies, $\pi_{t,e}$ depends on $e$ but not on $t$. We will also distinguish between several natural situations, namely
between strategies involving or not a sensing mechanism, and strategies $\pi_{t,e}$ depending or not depending on $t$, called hereditary and non-hereditary in what follows. We illustrate some results when both $\cT$ and $\cE$ are finite, and for 
\begin{enumerate}
\item{} non-hereditary strategies with no sensing mechanism, where $\pi_{t,e}\equiv p$, for some probability measure on $\cT$,
\item{} and for non-hereditary strategies using a sensing mechanism, of the form $\pi_{t,e}=p_e$, where $\bar p=(p_e)_{e\in\cE}$ is a family of probability measures on $\cT$. This second family  contains responsive switching strategies.
\end{enumerate}

We will describe the optimal set for non-hereditary strategies, that is, we will characterize the set of distributions
$(\pi_{t,e})$ yielding the fastest growth rate. We assume that the random environment $(\omega_n)_{n\ge 0}$ is a stationary
process taking values in $\cE$.
\subsubsection{ The non-hereditary with no sensing mechanism case\label{s.non-hereditary/no-sensing}} 

We first prove in Proposition \ref{prop:pr1} that
$$\lim_{n\to\infty}\frac{1}{n}\log \bbE_\omega[|Z_n|] =\gamma(p),$$
where the Lyapounov exponent $\gamma(p)$ is given by
$$\gamma(p) = \bbE [\log m_{p,\omega_0}],$$
and  $m_{p,e}$ is the first moment or average value of the mean distribution
$$\Upsilon_{p,e} = \int_{\cT} \Upsilon_{t,e}p({\rm d}t).$$
Let $\cP(\cT)$ be the set of probability measures on $\cT$ and $\gamma^* =\sup\{\gamma(p);\ p\in \cP (\cT)\}$ be the top Lyapounov exponent. We denote by
$\cP^*$ the subset of $\cP (\cT)$ containing the distributions $p$ maximizing
the Lyapounov exponent, i.e. $\gamma^*=\gamma(p)$.   We show in Proposition \ref{prop:pr3}
that $p\in \cP^*$ if and only if
$$\int_{\cE}\frac{m_{t,e}}{m_{p,e}}\nu_1({\rm d}e)\le 1,\ \forall t\in \cT.$$
where $\nu_1$ is the law of the stationary random environment. A strategy is called {\it pure}
if it is concentrated on a single $t\in\cT$, that is, takes the form $p=\delta_t$.  It then follows that
a pure strategy $p=\delta_t$ is optimal if and only if
$$\int_{\cE}\frac{m_{t',e}}{m_{t,e}}\nu_1({\rm d}e)\le 1,\ \forall t'\in \cT.$$
When $\cT =\{t_1,\cdots,t_q\}$ and  $\cE =\{e_1,\cdots,e_p\}$, with $q\ge p$, we also prove 
that there is a unique maximizer in $\cP^*$ which is supported by a set containing at most $p$ elements. 
\bigskip

As a further illustration, we consider the simplest case when $p=q=2$. 
\begin{enumerate}
\item{}  $p^\ast=\delta_{t_2}$ and $\gamma^\ast=\nu_1(e_1)\log(m_{t_2,e_1})+\nu_1(e_2)\log(m_{t_2,e_2})$ when $\frac{\nu_1(e_1)m_{t_1,e_1}}{m_{t_2,e_1}}+\frac{\nu_1(e_2)m_{t_1,e_2}}{m_{t_2,e_2}}\leq 1$,
\item{}    $p^\ast=\delta_{t_1}$ and $\gamma^\ast=\nu_1(e_1)\log(m_{t_1,e_1})+\nu_1(e_2)\log(m_{t_1,e_2})$
when $\frac{\nu_1(e_1)m_{t_2,e_1}}{m_{t_1,e_1}}+\frac{\nu_1(e_2)m_{t_2,e_2}}{m_{t_1,e_2}}\leq 1$.
\item{} Otherwise  $p^\ast$ is given by
$$p^\ast=\left(\frac{\nu_1(e_1)m_{t_2,e_2}}{m_{t_2,e_2}-m_{t_1,e_2}}+\frac{\nu_1(e_2)m_{t_2,e_1}}{m_{t_2,e_1}-m_{t_1,e_1}}\right)\delta_{t_1}+ \left(\frac{\nu_1(e_1)m_{t_1,e_2}}{m_{t_1,e_2}-m_{t_2,e_2}}+\frac{\nu_1(e_2)m_{t_1,e_1}}{m_{t_1,e_1}-m_{t_2,e_1}}\right)\delta_{t_2}$$
with
\begin{eqnarray*}
\gamma^\ast&=&\log|m_{t_1,e_1}m_{t_2,e_2}-m_{t_1,e_2}m_{t_2,e_1}|-\nu_1(e_1)\log|m_{t_2,e_2}-m_{t_1,e_2}|-\nu_1(e_2)\log|m_{t_2,e_1}-m_{t_1,e_1}|\\
&&+\nu_1(e_1)\log(\nu_1(e_1))+\nu_1(e_2)\log(\nu_1(e_2)).
\end{eqnarray*}
\end{enumerate}
\bigskip
As a numerical example, consider the case $\nu_1(e_1)=\nu_1(e_2)=1/2$ and 
$$m_{t_1,e_1}=1.5,\ m_{t_2,e_1}=0.6,\ m_{t_1,e_2}=0.6,\ m_{t_2,e_2}=1.5.$$
The two environments are equiprobable, with type $t_1$ better fitted to environment $e_1$, and a symmetry in the model. In this case, the Lyapounov exponent for the population without mutation is given by
$$\gamma(\delta_{t_1})=\gamma(\delta_{t_2})=0.5(\log(1.5)+\log(0.6))\approx -0.053.$$
This implies that both types are subcritical and that the corresponding homogeneous populations extinct almost surely.
The optimal strategy is then $p^\ast=(0.5,0.5)$ and the corresponding Lyapounov exponent is
$$\gamma^\ast=\log(0.5*(1.5+0.6))\approx 0.049.$$
With this optimal strategy, $\gamma^\ast>0$ and  the population explodes with positive probability almost surely. This is one of the simplest example when polymorphism is a necessary condition for survival (see section \ref{sec:ext1}).

\subsubsection{ The non-hereditary with  sensing mechanism case\label{s.non-hereditary/sensing}} 
Assume that $\pi_{t,e} =p_e$ for some strategy $\bar p = (p_e)_{e\in\cE}$. It turns out that the relevant piece of
environment is given by the pair process $\omega^{(2)}=((\omega_{n-1},\omega_n))_{n\ge 1}\in \cE^2$, of stationary measure
$\nu_2({\rm d}e_1,{\rm d}e_2)$. Consider the average distribution
$$\Upsilon_{\bar p,(e_1,e_2)}=\int_{\cT}\Upsilon_{t,e_2} p_{e_1}({\rm d}t),$$
of expected value
$$m_{p_{e_1},e_2} =\int_{\cT}m_{t,e_2}p_{e_1}({\rm d}t).$$
When the following integral exists
$$\gamma(\bar p)=\int_{\cE^2}\log (m_{p_{e_1},e_2})\nu_2({\rm d}e_1,{\rm d}e_2),$$
we prove in Proposition \ref{prop:pr5} that
$$\lim_{n\to\infty}n^{-1}\log \bbE_\omega[|Z_n|] =\gamma(\bar p)\quad \mbox{ almost\ surely.}$$

Let $\gamma^{**}$ be the optimal growth rate, that is, the supremum of $\gamma(\bar p)$ among all possible
families $\bar p$, and denote by $\cP^{**}$ the related set of optimal strategies. Let
$\nu_{e_1}({\rm d}e_2)$ be the conditional distribution of $\omega_2$ given that $\omega_1=e_1$.
We prove in Proposition \ref{prop:pr6b} that $\bar p$ is optimal if and only if
$$\int_{\cE}\frac{m_{t,e_2}}{m_{p_{e_1}},e_2}\nu_{e_1}({\rm d}e_2)\le 1,\ \forall t\in\cT,$$
$\nu_1({\rm d}e_1)$ a.s.
An interesting consequence is that there is no gain to be expected using a sensing mechanism when the random environment
has some independence property: If $\omega_1$ and $\omega_2$ are independent, then $\gamma^* =\gamma^{**}$.
The information one can gather on the present environmental state using sensors does not help when dealing with future events.

Next, we develop further the simplest case when $\cE=\{e_1,e_2\}$ and $\cT=\{t_1,t_2\}$. We suppose that the distribution for the environment $\nu$ is a Markov chain with transitions
\begin{eqnarray*}
\nu(\omega_{n+1}=e_1 |\omega_n=e_1)=1-q_1 &\quad& \nu(\omega_{n+1}=e_2 |\omega_n=e_1)=q_1 \\
\nu(\omega_{n+1}=e_1 |\omega_n=e_2)=q_2 &\quad& \nu(\omega_{n+1}=e_2 |\omega_n=e_2)=1-q_2 \\
\end{eqnarray*}
where $q_i\in (0,1)$ denotes the probability that the environment switches in the next step when it is currently in state $e_i$.
The sequence $\omega$ is then ergodic and stationary distribution 
$$\nu_1=\frac{q_2}{q_1+q_2}\delta_{e_1}+\frac{q_1}{q_1+q_2}\delta_{e_2}.$$
We suppose furthermore that the Markov chain is started from the steady state, i.e. $\omega_0$ has distribution $\nu_1$, so that the sequence $\omega$ is stationary. Then, the conditional distributions of $\omega_1$ given $\omega_0$ is precisely given by the transition of the Markov chain. The optimal strategy with sensing is then such that 
\begin{enumerate}
\item  $p^{\ast\ast}_{e_1}=\delta_{t_2}$ when  $\frac{(1-q_1)m_{t_1,e_1}}{m_{t_2,e_1}}+\frac{q_1m_{t_1,e_2}}{m_{t_2,e_2}}\leq 1$,
\item  $p^{\ast\ast}_{e_1}=\delta_{t_1}$ when  $\frac{(1-q_1)m_{t_2,e_1}}{m_{t_1,e_1}}+\frac{q_1m_{t_2,e_2}}{m_{t_1,e_2}}\leq 1$,
\item  and otherwise, 
$$p^{\ast\ast}_{e_1}=\left(\frac{(1-q_1)m_{t_2,e_2}}{m_{t_2,e_2}-m_{t_1,e_2}}+\frac{q_1m_{t_2,e_1}}{m_{t_2,e_1}-m_{t_1,e_1}}\right)\delta_{t_1}+ \left(\frac{(1-q_1)m_{t_1,e_2}}{m_{t_1,e_2}-m_{t_2,e_2}}+\frac{q_1m_{t_1,e_1}}{m_{t_1,e_1}-m_{t_2,e_1}}\right)\delta_{t_2},$$
\end{enumerate}
with similar equations  for $p^{\ast\ast}_{e_2}$. Rather than giving a general formula, we consider the numerical case when $q_1=q_2=q$ so that the stationary distribution is given by $\nu_1(e_1)=\nu_1(e_2)=1/2$ and we use the same values for the number of offspring as in Section \ref{s.non-hereditary/no-sensing}, that is we set
$$m_{t_1,e_1}=1.5,\ m_{t_2,e_1}=0.6,\ m_{t_1,e_2}=0.6,\ m_{t_2,e_2}=1.5.$$
Using the above results, we obtain that the optimal strategy $p^{\ast\ast}$ is given by:
\begin{eqnarray*}
&&p^{\ast\ast}_{e_1}=\delta_{t_1}\quad {\rm and}\quad p^{\ast\ast}_{e_2}=\delta_{t_2}\quad {\rm when} \quad q\leq \frac{2}{7}\\
&&p^{\ast\ast}_{e_1}=\delta_{t_2}\quad {\rm and}\quad p^{\ast\ast}_{e_2}=\delta_{t_1}\quad{\rm when} \quad q\geq \frac{5}{7}\\
&&p^{\ast\ast}_{e_1}=\frac{5-7q}{3}\delta_{t_1}+\frac{7q-2}{3}\delta_{t_2}\quad {\rm and}\quad p^{\ast\ast}_{e_2}=\frac{7q-2}{3}\delta_{t_1}+\frac{5-7q}{3}\delta_{t_2}\quad {\rm if}\quad \frac{2}{7}\leq q\leq \frac{5}{7}.
\end{eqnarray*}
\medskip

We see three different environmental regimes, corresponding to low, intermediate and high frequency switching rates. When the environment fluctuates slowly, 
with $q\leq \frac{2}{7}$, the optimal strategy is pure and corresponds to what \cite{KuLei} called responsive switching. In the intermediate regime, the optimal strategy is a mixture of pure strategies. In the high frequency regime  where environment fluctuates quickly ($q\ge 5/7$), the optimal strategy is a pure one, where organisms being in the favorable state are pushed to the unfavorable state.
\medskip

We deduce the optimal growth rate
$$\gamma^{\ast\ast}=\left\{
\begin{array}{ll}\log\frac{3}{2}-q\log\frac{5}{2}&{\rm if}\ q\leq \frac{2}{7}\\
q\log q +(1-q)\log(1-q)+\log\frac{21}{10}  &{\rm if}\  \frac{2}{7}\leq q\leq \frac{5}{7}\\ 
\log\frac{3}{5}-q\log\frac{5}{2}&{\rm if}\ q\geq \frac{5}{7}
\end{array}
\right.. $$
When $q=1/2$, the sequence $\omega$ is an independent sequence, so that $\gamma^{\ast\ast}(1/2)=\gamma^\ast(1/2)\approx 0.049$, the optimal strategy is a strategy without sensing. When $q\approx 1$ (resp. $q\approx 0$), the environment in the next step is very likely to stay the same (resp. to switch), so that we can determine with high probability which type will be fitted in the next generation. We observe indeed that
$$\lim_{q\to 0} \gamma^{\ast\ast}(q)=\lim_{q\to 1} \gamma^{\ast\ast}(q)=\log\frac{3}{2}\approx 0.176.$$

\section{Definition of the Multitype Branching Process in Random Environment and of optimal strategies}

First, we recall that  $\omega_n$ represents  the environment at time $n$, $Z_n$ the trait-structured population, $Z_n^t$ the number of individuals with trait $t$ and $|Z_n|$ the total number of individuals at time $n$. Moreover  $\xi_{i,n}^{t,\omega_n}$ gives the size of the offspring of the $i$-th individual with trait $t$ at time $n$ in environment $\omega_n$, and $\tau_{i,j,n}^{t,\omega_n}$ gives  the trait of the $j$-th descendant of this individual. \\

Let us give now the formal definition of the process. \\
Let $(\Omega,\cF,\bbP)$ be a probability space and $\cT$ be a metric space and $\cE$ be a Polish space. \\
For each pair $(t,e)\in\cT\times\cE$, let $\Upsilon_{t,e}$ and $\pi_{t,e}$ be distributions on $\bbN$ and $\cT$ respectively and  suppose that $\Upsilon_{t,e}$ has a finite first moment denoted by $m_{t,e}$. \\
Let $\omega=(\omega_n)_{n\geq 0}$ be a $\cE$-valued stationary ergodic random process with distribution $\nu$ on $\cE^{\bbN}$.  \\
Denoting by $\bbN^{\cT}$  the set of almost null $\bbN$-valued functions on $\cT$ and  $(1_t)_{t\in \cT}$  the canonical basis,  the stochastic process $Z_n=(Z_n^{t} : t \in\cT)$ with values in $\bbN^{\cT}$ can be defined as:
\begin{equation}\label{model}
  \left\{
      \begin{aligned}
       &Z_0=N_0,\\
     &Z_{n+1}=\sum_{t\in\cT}\sum_{i=1}^{Z_n^{t}} \sum_{j=1}^{\xi_{i,n}^{t,\omega_n}} 1_{\tau_{i,j,n}^{t,\omega_n}},\quad n\geq 0,\\
      \end{aligned}
  \right.
\end{equation}
where
\begin{itemize}
\item $N_0$ is a $\bbN^{\cT}$-valued random variable giving the population at time $0$,
\item  $\{\xi_{i,n};i\geq 1,n\geq 0\}$ and $\{\tau_{i,j,n};i\geq 1,j\geq 1,n\geq 0\}$ are infinite arrays of iid random variables with values in $\bbN^{\cT\times\cE}$ and $\cT^{\cT\times\cE}$ respectively and common distribution 
$$\Upsilon=\underset{(t,e)\in \cT\times\cE}{\otimes}\Upsilon_{t,e}, \quad  \pi=\underset{(t,e)\in \cT\times\cE}{\otimes} \pi_{t,e}$$ respectively,
\item $N_0$, $\{\xi_{i,n};i\geq 1,n\geq 0\}$, $\{\tau_{i,j,n};i\geq 1,j\geq 1,n\geq 0\}$ and $\omega=(\omega_n)_{n\geq 0}$ are independent.
\end{itemize}
$\newline$

 Thus $\Upsilon$ corresponds to the offspring distribution  and $\pi$ to the distribution of trait of the offspring. Both distributions a priori depend on the trait of the parent and on the environment. \\
We  put  $|Z_n|=\sum_{t\in\cT}Z_n^t,\quad n\geq 0$ (an empty sum is equal to zero). \\

Our attention will be focused on the role of the trait distribution $\pi$ and how it affects the evolution of the population. The offspring distribution $\Upsilon$ will hence be considered as given and fixed, whereas the trait distribution $\pi$ will be considered  as a parameter. The intuitive idea is that i) the local fitness of an individual is determined by its trait and its environment and corresponds to its mean number of descendants,  ii) the intergenerational variability of the traits has to be tuned so as to maximize the global fitness of the whole population.

The trait distribution $\pi$ can be seen as the strategy used by the population to maximize its growth. We measure the performance of the strategy $\pi$ by the long term growth of the corresponding population $(Z_n)_{n\geq 0}$: let 
$$\gamma(\pi)=\lim_{n\to\infty} n^{-1}\log \bbE_{\omega,N_0}\left[|Z_n|\right].$$
Here $\bbE_{\omega,N_0}$ denotes conditional expectation given the environment and the initial population.
As will be precised in the sequel, fairly general assumptions ensures that this limit exists and does not depend on the environment and initial population. Note that it is important here to consider the quenched model (i.e. conditionally on the environment) and note the averaged one (i.e. averaged on all possible environments): the criteria $\gamma(\pi)$ measures if the population behaves well in a typical environment. An averaged criteria would be biased by unlikely environments where the population grows unusually faster.

Several mechanisms can induce the intergenerational variability of the traits and leads to different assumption on the trait distribution $\pi$. \\

\section{The non-hereditary case with no sensing mechanism}
According to the approach set out in the previous paragraph, we begin by studying the simplest case of non-hereditary traits in absence of sensing mechanism, meaning that the trait distribution of the offspring depends neither on the trait of the parent nor on the environment. We thus suppose that $Z_n$ evolves according to model \eqref{model} with $\pi_{t,e}\equiv p$ for some distribution $p$ on $\cT$, and let $\pi=\pi(p)$ be the corresponding product trait distribution. 

\subsection{Reduction to a simple BPRE}
In some sense, the non-hereditary assumption makes the structure of the population trivial: the trait distribution is given by $p$ whatever the past evolution of the process was. In mathematical terms, we take advantage of some stochastic independence:

\begin{lem}\label{lem1}
For any $n\geq 1$, the population structure is conditionally independent of the past population process given the size of the population, i.e.
$$Z_{n}\Big ||Z_{n}| \quad \coprod \quad (Z_{0},\cdots,Z_{n-1})\Big ||Z_{n}|.$$ 
\end{lem}
\noindent 
{\it Proof.}\ It is easily seen from the assumptions on the model \eqref{model} and from the non-hereditary assumption $\pi=\pi(p)$ that the distribution of $Z_{n}$  given $(Z_0,\cdots Z_{n-1})$ and $|Z_{n}|$ is equal to the distribution of $\sum_{1\leq i\leq |Z_{n}|} 1_{\tau_{i}}$ with $\tau_i$ an iid sequence with distribution $p$. This distribution does not depend on $(Z_0,\cdots,Z_{n-1})$, this proves the conditional independence. \CQFD
\\ \ \\

It is worth noting that the result also holds for the quenched model (i.e. conditionally on the environment $\omega$). The Lemma implies that the distribution of the population process $(Z_n)_{n\geq 1}$ is easily recovered from the size process $(|Z_n|)_{n\geq 1}$. This latter process turns out to be a simple BPRE and this allows us to compute the performance $\gamma(p)$ of the strategy $p$.

\begin{pr}\label{prop:pr1}
The size process $(|Z_n|)_{n\geq 1}$ is a simple BPRE with offspring distribution in environment $e$ given by
$$ \Upsilon_{p,e}= \int_{\cT} \Upsilon_{t,e}p(dt),\quad e\in\cE.$$
Conditionaly to $\omega$, the expected population size at time $n$ is 
$$\bbE_\omega[|Z_n|]=\bbE_{\omega_0}[|Z_1|]\prod_{k=1}^{n-1} m_{p,\omega_k}$$
with $ m_{p,e}$ the first moment of $\Upsilon_{p,e}$. If $\gamma(p)=\bbE\left[\log m_{p,\omega_0}\right]$ exists, then 
$$\lim_{n\to\infty}n^{-1}\log \bbE_\omega[|Z_n|] =\gamma(p)\quad \mbox{ almost\ surely.}$$
\end{pr}
\noindent
{\it Proof.}\ According to Lemma \ref{lem1}, given $(|Z_1|,\cdots,|Z_n|)$, the population $Z_n$ has the same distribution as   
$\sum_{1\leq i\leq |Z_{n}|} 1_{\tau_{i}}$ with $\tau_i$ an iid sequence with distribution $p$. Intuitively, the $i$-th individual  has type $\tau_i$ chosen randomly on $\cT$ with distribution $p$. The size of the next generation $|Z_{n+1}|$ is then $\sum_{1\leq i\leq |Z_n|}\xi_{n,i}^{\tau_i,\omega_n}$ where $\xi_{n,i}^{\tau_i,\omega_n}$ is the offspring of the $i$-th individual of type $\tau_i$ in environment $e_n$. From this two step procedure i) random choice of the type $t$, ii) reproduction with random offspring in environment $e$, we obtain the effective offspring distribution $\Upsilon_{p,e}$ in environment $e$: it is the mixture of the offspring distributions $\Upsilon_{t,e}$, with mixing distribution $p$.
This proves the branching property for $(|Z_n|)_{n\geq 1}$. The other properties follow: the branching property implies the recursive formula
$$ \bbE_\omega[|Z_{n+1}|]=\bbE_{(\omega_0,\cdots,\omega_{n-1})}[|Z_n|] m_{p,\omega_n}$$ 
and we get the formula for $\bbE_\omega[|Z_{n}|]$ follows. Taking the logarithm, we have
$$n^{-1}\log \bbE_\omega[|Z_n|]=n^{-1}\log \bbE_{\omega_0}[|Z_1|]+n^{-1}\sum_{k=1}^{n-1} \log m_{p,\omega_k}$$
which converges to $\gamma(p)$ almost surely according to the ergodic theorem and the integrability assumptions.
\CQFD
\\ \ \\

For the sake of simplicity, we suppose in what follows that the initial population consists in a single individual with random trait with distribution $p$. The whole process $(|Z_n|)_{n\geq 0}$ is then a simple BPRE with offspring distribution $\Upsilon_p$ and initial condition $|Z_n|=1$. Simple branching processes in random environment have been introduced by Smith-Wilkinson \cite{SmWi} and Athreya-Karlin \cite{AtKa1,AtKa2} and have been studied rather intensively since then \cite{Tan,GKV, AGKV}. We recall here some important results concerning the asymptotic behaviour of such processes (see e.g. the classification Theorem in \cite{Tan}): it states conditions under which either the population becomes extinct or explodes at geometric rate. 

To avoid the trivial case of a constant population, we suppose that $\Upsilon_{p,\omega_0}$ is not  almost surely equal to $\delta_1$. We suppose also that $\gamma(p)=\bbE\left[ \log m_{p,\omega_0}\right]$ exists and is finite.
We say  that extinction occurs if $Z_n\to 0$ as $n\to\infty$ (in which case the sequence vanishes eventually), and that the population survives otherwise. Let $q(\omega)=\bbP_\omega\left[Z_n\to 0\right]$ be the probability of extinction given the environment.
\begin{theo}\label{theo1} (classification theorem) 
\begin{enumerate}
\item In both subcritical case $\gamma(p)<0$ and critical case $\gamma(p)=0$, the population becomes extinct almost surely, i.e. $\bbP\left[ q(\omega)=1\right]=1$
\item In the supercritical case $\gamma(p)>0$, if furthermore $\bbE\left[-\log\left(1-\Upsilon_{p,\omega_0}(\{0\})\right)\right]<\infty$, then the population can survive with positive probability in almost every environment, i.e. $\bbP\left[ q(\omega)<1\right]=1.$\\
Furthermore, conditionally on nonextinction, the population explodes at a geometric rate:
$$\lim_{n\to\infty} n^{-1}\log Z_n=\gamma(p)\quad \mbox{almost\ surely\ on}\quad \{\forall n\geq 0,\ Z_n>0\}.$$
\end{enumerate}
\end{theo}

\subsection{Optimal strategies}
We now focus on optimal strategies: i.e. what choice of the distribution $p$ allows for the fastest growth of the population ? 
As discussed in Section $2$, the performance of the strategy $\pi=\pi(p)$ is measured by the Lyapounov exponent $\gamma(p)$. 
In the non-hereditary case with no sensing, an explicit formula for $\gamma(p)$ has been derived in Proposition \ref{prop:pr1}.
The question naturally arise to determine the supremum of $\gamma(p)$ for $p$ varying in $\cP(\cT)$ the space of distribution on $\cT$ and the set of optimal strategies $p^\ast$ in the case when this supremum is reached. Let $\gamma^\ast=\sup\{\gamma(p); \ p\in\cP(\cT)\}$ be the optimal growth rate and $\cP^\ast=\left\{p^\ast\in\cP(\cT);\ \gamma(p^\ast)=\gamma^\ast\right\}$ the set of optimal strategies. A strategy $p$ is called pure if $p=\delta_t$ for some $t\in\cT$ (i.e. all individual in the population have the same trait $t$) or mixed otherwise. One ask further whether optimal strategies are pure or mixed.

We now precise conditions under which existence of optimal strategies are ensured. We suppose that:
\begin{itemize}
\item[$({\bf C}_1)$]\quad there is some $M>0$ such that $ m_{t,e}\leq M$ for all $(t,e)\in\cT\times\cE$,
\item[$({\bf C}_2)$]\quad for all $e\in\cE$, the application $t\mapsto m_{t,e}$ is continuous on $\cT$,
\item[$({\bf C}_3)$]\quad for all $e\in\cE$ and $\varepsilon>0$, there is a compact set $K\subset\cT$ such that $m_{t,e}\leq \varepsilon$ for all $t\in\cT\setminus K$
\end{itemize}
Condition $({\bf C}_1)$ is rather relevant from the biological point of view since an individual could hardly have arbitrary high number of offsprings in a fixed amount of time. Conditions $({\bf C}_2)$ and $({\bf C}_3)$ are related with the topology of $\cT$: individuals with close traits are supposed to have approximately the same behaviour, and traits close to infinity are supposed to have a poor fitness.

\begin{pr}\label{prop:pr2}
Under conditions $({\bf C}_1)-({\bf C}_3)$, optimal strategies exists and form a closed convex set i.e. the set 
$\cP^\ast$ is nonempty closed and convex in $\cP(\cT)$ endowed with the topology of weak convergence.\\ 
If furthermore the family  $\cM=\{e\mapsto m_{t,e}; t\in \cT\}\subset L^\infty(\cE,\nu_1)$ is linearly independent, then the optimal strategy is unique.
\end{pr}

\noindent
{\it Proof.} First, we prove that the map $\gamma:\cP(\cT)\to \bbR$ defined by
\begin{equation}\label{eq:propgamma}
\gamma(p)=\int_{\cE}\log\left(m_{p,e}\right)\nu_1(de)\quad \mbox{with}\quad m_{p,e}=\int_{\cT} m_{t,e}p(dt)
\end{equation}
is concave and upper semi-continuous on $\cP(\cT)$ with respect to the topology of weak convergence. 
Let $p_n$ converging weakly to $p$ in $\cP(\cT)$.  Using conditions $({\bf C}_1)$ and $({\bf C}_2)$, we see that $m_{p_n,e}\to m_{p,e}$ for all $e\in \cE$. Define $f_n(e)=\log(M)-\log\left(m_{p_n,e}\right)$ be non-negative measurable functions on $\cE$ and apply Fatou's Lemma: we get
$$\int_{\cE} \liminf f_n(e)\nu_1(de)\leq \liminf \int_{\cE} f_n(e)\nu_1(de).$$
Equivalently,
$$\limsup \int_{\cE}\log\left(m_{p_n,e}\right)\nu_1(de)\leq  \int_{\cE}\log\left(m_{p,e}\right)\nu_1(de).$$
Hence $\limsup \gamma(p_n)\leq \gamma(p)$ and this states the upper semicontinuity of the application $\gamma$. The concavity of the application $\gamma$ on $\cP(\cT)$ is a direct consequence from the concavity of the logarithm and the linearity of integration: first for $p_1,p_2\in\cP(\cT), \lambda\in [0,1]$ and $e\in \cE$, we have
$$\log\left(m_{\lambda p_1+(1-\lambda)p_2,e}\right)=\log\left(\lambda m_{p_1,e}+(1-\lambda)m_{p_2,e}\right)\geq 
\lambda \log\left( m_{p_1,e}\right)+(1-\lambda)\log\left(m_{p_2,e}\right).$$
This implies $\gamma(\lambda p_1+(1-\lambda)p_2)\geq \lambda \gamma(p_1)+(1-\lambda)\gamma(p_2)$.

We then prove that $\cP^\ast$ is nonempty closed and convex. The closeness and convexity properties are straightforward since
$\cP^\ast$ can be seen as the level set $\{p\in\cP(\cT); \gamma(p)\geq p^\ast\}$ of the concave uppersemicontinuous application
$\gamma$. It remains to check non-emptiness and we use compactness arguments. In the case when $\cT$ is a compact space, then
$\cP(\cT)$ is also compact with respect to the weak topology and the upper semicontinuous map $\gamma$ reaches its maximum on
$\cP(\cT)$ so that $\cP^\ast$ is non-empty. In the case when $\cT$ is non-compact, we consider its compactification
$\widehat\cT=\cT\cup\{\infty\}$. We extend the definition of $m$ by $m_{\infty,e}=0$ for all $e\in\cE$. Condition $({\bf C}_3)$
ensures that this extension is continuous on $\hat \cT$. Then the map $\gamma$ extends on the compact space $\cP(\hat\cT)$, is
uppersemicontinuous and hence reaches its maximum at some point $p^\ast\in \cP(\hat\cT)$. It remains to check that
$p^\ast(\{\infty\})=0$ so that $p^\ast$ can be seen as an element of $\cP(\cT)$. It is straightforward since
$p^\ast(\{\infty\})>0$ would imply $\gamma(p^\ast)=-\infty$.

At last, we prove the uniqueness in the case when the family $\cM$ is linearly independent. Let $p_1^\ast$ and $p_2^\ast$ be two points where $\gamma$ reaches its maximum. Using the strict concavity of the logarithm, we see that necessarily $m_{p^\ast_1,e}=m_{p^\ast_2,e}$ $\nu_1(de)$-almost everywhere. Using the linear independence, this in turn implies $p^\ast_1=p^\ast_2$.  \CQFD\\

The following characterization of optimal strategies can be useful: 
\begin{pr}\label{prop:pr3}
A strategy $p$ is optimal if and only if 
$$\int_{\cE} \frac{m_{t,e}}{m_{p,e}}\nu_1(de)\leq 1,\quad \forall t\in \cT.$$
\end{pr}
\noindent
{\it Proof.} The strategy $p$ is optimal if and only if 
$$\gamma((1-\varepsilon)p+\varepsilon p')-\gamma(p)\leq 0\quad  {\rm for\ all\ } p' {\rm\ and\  } 0<\varepsilon<1.$$
Using concavity and differentiability, this is equivalent to 
$$\frac{{\rm d}}{{\rm d}\varepsilon}\left[\gamma(p+\varepsilon(p'-p)\right]_{|\varepsilon=0}\leq 0,$$
which can be rewritten as  
$$\frac{{\rm d}}{{\rm d}\varepsilon}\left[\int_{\cE}\log\left(m_{p,e}+\varepsilon (m_{p',e}-m_{p,e})\right)\nu_1(de)\right]_{|\varepsilon=0}=\int_{\cE}\frac{m_{p',e}-m_{p,e}}{m_{p,e}}\nu_1(de)\leq 0.$$
Thus a necessary and sufficient condition is 
$$\int_{\cE} \frac{m_{p',e}}{m_{p,e}}\nu_1(de)\leq 1\quad ,\quad \forall p\in \cP(\cT).$$
It is easily seen that it is equivalent to test the condition for $p'=\delta_t,t\in\cT$ and this proves the announced result.\CQFD\\

As a direct application of Proposition \ref{prop:pr3}, we can answer the question wether there is a pure optimal strategy: 
\begin{cor}\label{cor1} The pure strategy $p=\delta_{t}$ is optimal if and only if
$$\int_{\cE} \frac{m_{t^{'},e}}{m_{t,e}}\nu_1(de)\leq 1\quad ,\quad \forall t^{'}\in\cT.$$
\end{cor}
This gives us a simple criterion for the existence of optimal pure strategies.

\subsection{Extinction}\label{sec:ext1}
According to Theorem \ref{theo1}, on the one hand the population becomes extinct in the case when $\gamma(p)\leq 0$ and on the other hand (under additional technical conditions) the population survives and explodes with positive probability when $\gamma(p)>0$. For the population's survival, strategies $p$ with $\gamma(p)>0$ are of highest importance. This leads to the definition of the set $\cS$ of strategies that allow the population to survive:
$$\cS=\{p\in\cP(\cT);\ \gamma(p)>0\}.$$
The following property holds:
\begin{pr}\label{prop:pr4}
If $\gamma^\ast\leq 0$, the set $\cS$ is empty. If $\gamma^\ast> 0$, $\cS $ is a nonempty convex set containing $\cP^\ast$.
\end{pr}
\noindent
{\it Proof.} From the definition of $\gamma^\ast$, $\gamma(p)\leq \gamma^\ast$ for all strategy  $p$. Hence $\cS$ is empty if $\gamma^\ast\leq 0$. If $\gamma^\ast>0$, then every optimal strategy $p^\ast\in\cP^\ast$ satisfies $\gamma(p^\ast)=\gamma^\ast>0$ and hence $\cP^\ast\subset\cS$. The map $p\mapsto \gamma(p)$ is concave function so that the level set $\cS=\{p\in\cP(\cT);\ \gamma(p)>0\}$ is convex.
\CQFD\\

In some cases, a striking phenomenon may happen:  no pure strategy can allow for survival, i.e. every  homogeneous population with a single trait $t$ suffers from extinction; but some mixed strategies may prevent from extinction, i.e. some  polymorphic populations may survive forever. In this case, we should say that polymorphism is a necessary condition for survival. This phenomenon occurs when $\gamma(\delta_t)\leq 0$ for all $t\in\cT$ whereas $\gamma^\ast>0$.

The intuitive idea is the following: extinction occurs when the environment is bad for almost all the individuals in the population and hence a diversification of the traits in the population should imply a smaller number of environments that are bad for almost all individuals. This can be seen as a consequence of the concavity property: suppose that for all trait $t$, $\gamma(\delta_t)\equiv \gamma$, i.e.  homogeneous population have the same growth rate no matter the trait $t$. Then for any strategy $p$, $\gamma(p)\geq \int_{\cT}\gamma(t)p(dt)=\gamma$, i.e. any polymorphic population has a better growth rate than any homogeneous population. See the examples below for further illustration of this phenomenon.

\subsection{Example: finite-dimensional case}\label{sec:fidi}
We consider the case when the vectorial space spanned  by the family $\cM=\{e\mapsto m_{t,e}, t\in\cT\}$ in $L^\infty(\cE,\nu_1)$ has finite dimension denoted by $d$. This occurs in particular as soon as: \\
- either the environment space is finite: $\cE=\{e_1,\cdots,e_p\}$ and $\nu_1$ is a discrete measure such that $\nu_1(e_i)>0$ for all $i\in\{1,\cdots,p\}$. In this case $L^\infty(\cE,\nu_1)$ is of dimension $p$ and $d\leq p$;\\
-  or the trait is finite: $\cT=\{t_1,\cdots,t_q\}$. Then $d\leq q$ and conditions $({\bf C}_1)-({\bf C}_3)$ are automatically satisfied.\\
For convenience, we require furthermore that the following conditions holds: 
\begin{itemize}
\item[$({\bf H})$] for any pairwise distinct $t_1,\cdots,t_d\in\cT$, the family of functions $\{e\mapsto m_{t_i,e};1\leq i\leq d\}$ is linearly independent in $L^\infty(\cE,\nu_1)$.
\end{itemize}
This property will be convenient because it ensures the uniqueness of the optimal strategy. More precisely,
\begin{pr}\label{prop:fidi1}
Suppose that ${\rm span}\cM$ has dimension $d$ and that conditions $({\bf C_1})-({\bf C_3})$ and $({\bf H})$ hold. Then:
\begin{itemize}
\item there exists a unique optimal strategy $p^\ast\in\cP(\cT)$,
\item  $p^\ast$ is a discrete probability on $\cT$ supported by at most $d$ different types.
\end{itemize}
\end{pr}
\noindent
{\it Proof.} Existence of an optimal strategy is a consequence of Proposition $\ref{prop:pr2}$. However it is worth noting that in this finite-dimensional case, it is a consequence of standard analysis. We denote by ${\rm Conv}(\cM)$ the closed convex hull of $\cM$ in $L^\infty(\cE,\nu_1)$ which can be seen as a closed convex set in a $d$-dimensional vector space. Note that ${\rm Conv}(\cM)$ is equal to the set of functions $\{m_p:e\mapsto m_{p,e};\ p\in\cP(\cT)\}$ where $m_{p,e}$ is defined in equation \eqref{eq:propgamma}. Introduce the application $\Theta:{\rm Conv}(\cM)\to [-\infty,\infty)$ defined by $\Theta(m)=\int_{\cE}\log(m(e))\nu_1(de)$ (with the convention that $\log(u)=-\infty$ if $u\leq 0$). With these notations, we easily see that $\gamma(p)=\Theta(m_p)$ and maximizing $\gamma$ on $\cP(\cT)$ is equivalent to maximizing $\Theta$ on ${\rm Conv}(\cM)$. Next we observe that the function $\Theta$ is uppersemicontinuous and strictly concave, so that it reaches its maximum at a unique point $m^\ast$ on the compact set  ${\rm conv}(\cM)$. Furthermore the gradient of $\Theta$ never vanishes so that the extremum must be reached on a boundary point of ${\rm Conv}(\cM)$ and  $m^\ast\in\partial {\rm Conv}(\cM)$. The boundary of ${\rm Conv }(\cM)$ consists in extremal types (type $t$ such that $m_t:e\mapsto m_{t,e}$ is an extremal point of the convex set ${\rm Conv}(\cM)$) and of $k$-dimensional faces determined by the convex hull of $k+1$ extremal types, with $0\leq k\leq n$. If $m^\ast$ belongs to such a $k$-dimensional face, then $m^\ast=m_{p^\ast}$ where $p^\ast$ is a discrete probability on $\cT$ supported by the $k+1$ corresponding extremal types, with $k+1\leq d$. The uniqueness  property comes then from assumption $({\bf H})$. Indeed, under this assumption, a $k$ dimensional face determined by $k+1$ extremal types does not contain any other point of $\cM$ so that the decomposition in barycentric coordinates are unique. \CQFD 
\\ \ \\
It is worth noting that the application $\Theta$ has always a unique minimizer $m^\ast$ in the above proof. Condition $({\bf H})$ ensures that there is a unique strategy $p^\ast$ such that $m_{p^\ast}=m^\ast$. In absence of this condition, there might be several mixing distribution $p$ such that $m_p=m^\ast$ and then several optimal strategies.
\\ \ \\
An appealing particular case is the case when both $\cE=\{e_1,\cdots, e_p\}$ and $\cT=\{t_1,\cdots,t_q\}$ are finite with $q\geq p$. In this case $L^\infty(\cE,\nu_1)$ is of dimension $d=p$, and a generic configuration will always satisfy assumption $({\bf H})$ (in the sense that an exact linear relation between the types is very unlikely from a biological point of view). Then, according to Proposition \ref{prop:fidi1}, there exists a unique optimal strategy mixing at most $p$ different types. This shows that the number of types involved in the optimal strategy is less than the number of different environments. This is reminiscent from a rule in ecology stating that the number of species in an ecosystem is bounded above by the numbers of different niche in the sense that  two species cannot occupy the same niche for a long time (competitive exclusion principle). From a practical point of view, in this finite settings, the optimal strategy can be computed using standard numerical convex optimization (see \cite{BG} for instance). For more illustrations in this setting, see Section \ref{s.non-hereditary/no-sensing}.
\\ \ \\

\subsection{Example: Gaussian distributions}
This example is due to Haccou and Iwasa \cite{HacIwa}: the environment space is the set of real numbers $\cE=\bbR$,  the environment $\omega$ is supposed to be a Gaussian stationary ergodic sequence with stationary distribution $\nu_1=\cN(\mu,\sigma_2^2)$ the Gaussian distribution with mean $\mu\in\bbR$ and variance $\sigma^2_2>0$. The dynamic for the environment is irrelevant in the no-sensing case. The trait space is the set of real numbers $\cT=\bbR$ and the mean offspring number of an individual of trait $t$ in environment $e$ has the gaussian form $$m_{t,e}=\frac{C}{\sqrt{2\pi\sigma_1^2}}\exp\left(-\frac{(t-e)^2}{2\sigma_1^2}\right)$$
for some parameters $C>0$ and $\sigma_1^2>0$. In environment $e$, individuals with trait value $t=e$ are the best fitted.
Note that conditions $({\bf C}_1)-({\bf C}_3)$ are fulfilled so that Proposition \ref{prop:pr2} holds: optimal strategies exist. Furthermore, the family $\cM=\{e\mapsto m_{t,e};t\in\cT\}$ is linearly independent so that unicity hold: there is a unique optimal strategy $p^\ast$ depending a priori on $\mu,\sigma_1,\sigma_2$ and $C$.

First we easily compute the fitness of a pure strategy
$$\gamma(\delta_t)=\int_{\cE}\log(m_{t,e})\nu_1(de)= \log C-\frac{1}{2}\log(2\pi\sigma_1^2)-\frac{(\mu-t)^2+\sigma_2^2}{2\sigma_1^2}$$
and the optimal pure strategy is equal to $\delta_\mu$.
Then according to Corollary \ref{cor1}, we test if this pure strategy is optimal in the set of all mixed strategies. For this we compute
$$\int_{\cE}\frac{m_{t,e}}{m_{\mu,e}}\nu_1(de)= \int_{\cE} \exp\left(\frac{\mu^2-t^2+2e(t-\mu)}{2\sigma_1^2} \right)\nu_1(de)=\exp\left(\frac{(t-\mu)^2(\sigma_2^2-\sigma_1^2)}{2\sigma_1^4} \right)$$
and check that it is less than one when $\sigma_2\leq \sigma_1$. Hence if the fluctuations of the environment are small ($\sigma_2\leq\sigma_1$), then the pure strategy $\delta_\mu$ is optimal. 

Gaussian strategies $p=\cN(\mu_p,\sigma^2_p)$ yield easy computations: recalling that a Gaussian mixture of Gaussian distributions is a again Gaussian, we compute
$$m_{p,e}=\int_{\cT}m_{t,e}p(dt)=\frac{C}{\sqrt{2\pi(\sigma_1^2+\sigma_p^2)}}\exp\left(-\frac{(\mu_p-e)^2}{2(\sigma_1^2+\sigma_p^2)}\right)$$
and then
$$\gamma(p)=\int_{\cE}\log(m_{t,e})\nu_1(de)= \log C-\frac{1}{2}\log(2\pi(\sigma_1^2+\sigma_p^2))-\frac{(\mu-\mu_p)^2+\sigma_2^2}{2(\sigma_1^2+\sigma_p^2)}.
$$
We find that $\gamma(p)$ is maximal (among Gaussian strategies) for $\mu_p=\mu$ and 
$$\sigma_p^2=\left\{\begin{array}{ll} 0 & {\rm\ if\ } \sigma_2^2\leq \sigma_1^2 \\ \sigma_2^2-\sigma_1^2 & {\rm\ if\ }  \sigma_2^2> \sigma_1^2 \end{array}\right..$$
We have seen that this is the optimal strategy in the case when $\sigma_2^2\leq \sigma_1^2$. This is still the case when 
$\sigma_2^2> \sigma_1^2$: we compute indeed for $t\in\cT$ and $p^\ast=\cN(\mu,\sigma_2^2-\sigma_1^2)$
$$ \int_{\cE}\frac{m_{t,e}}{m_{p^\ast,e}}\nu_1(de)= \int_{\cE} \frac{1}{C}m_{t,e}de=1$$
and this characterizes the optimal strategy according to Proposition \ref{prop:pr3}.

We have proven so far that:
\begin{pr}\label{pr:HaccIsa}
 the optimal strategy is
\begin{equation}\label{eq:optstrat}
p^\ast=\left\{\begin{array}{ll} \delta_\mu & {\rm\ if\ } \sigma_2^2\leq \sigma_1^2 \\ \cN(\mu,\sigma_2^2-\sigma_1^2) & {\rm\ if\ }  \sigma_2^2\geq \sigma_1^2 \end{array}\right.,
\end{equation}
and the optimal growth rate
\begin{equation}\label{eq:optgrow}
\gamma^\ast= \left\{\begin{array}{ll} 
\log C-\frac{1}{2}\log(2\pi\sigma_1^2)-\frac{\sigma_2^2}{2\sigma_1^2} & {\rm\ if\ } \sigma_2^2\leq \sigma_1^2 \\ 
\log C-\frac{1}{2}\log(2\pi\sigma_2^2)-\frac{1}{2} & {\rm\ if\ }  \sigma_2^2\geq \sigma_1^2 
\end{array}\right.. 
\end{equation}
\end{pr}
An interesting quantity is the relative gain of the best mixed strategy over the best pure strategy: it gives an indication of the strength of the selection pressure on mixed as opposed to pure strategies. From the above computation, denoting by  $\chi=\sigma_2^2/\sigma_1^2$, we get 
$$\gamma(p^\ast)-\gamma(\delta_\mu)=\left\{\begin{array}{ll} 
0 & {\rm\ if\ } \chi\leq 1 \\ 
\frac{1}{2}(\chi-1-\log\chi) & {\rm\ if\ } \chi\geq 1  
\end{array}\right. $$
This is the log-ratio of the expected long-term growth rates of individuals playing the different types of strategies.  It is non-decreasing with respect to $\chi$: when the environmental variance $\sigma_2^2$ is large compared to $\sigma_1^2$, there is a strong advantage in playing a mixed strategy. 

Next we illustrate the phenomenon discussed in section \ref{sec:ext1} when polymorphism is a necessary condition for survival. This happens when the optimal pure strategy $\delta_\mu$ leads to almost sure extinction of the population ($\gamma(\delta_\mu)\leq 0$) whereas the optimal mixed strategy $p^\ast$ allows the population to survive ($\gamma^\ast>0$).
This phenomenon occurs when $\sigma_2>\sigma_1$ for every $C$ belonging to the following non empty interval:
$$\frac{1}{2}\log(2\pi\sigma_2^2)+\frac{1}{2}< \log C\leq \frac{1}{2}\log(2\pi\sigma_1^2)+\frac{\sigma_2^2}{2\sigma_1^2}.$$

\section{The non-hereditary case with sensing mechanism}
We now study the case of non-hereditary traits  when some sensing-mechanism is available: we mean that the trait distribution of the offspring does not depends on the trait of the parent, but does depend on the environment because the individuals get some information about the environment they evolve in and are able to adapt suitably the trait distribution of their offspring. We thus suppose in this section that $Z_n$ evolves according to model \eqref{model} with   $\pi_{t,e}\equiv p_e$ for some family $\bar p=(p_e)_{e\in\cE}$ of distributions on $\cT$ and let $\pi=\pi(\bar p)$ be the corresponding product distribution. 

\subsection{Reduction to a simple BPRE}
In a similar way as in the non-hereditary case with no sensing, the non-hereditary assumption makes the structure of the population very simple: the trait distribution at time $n$ is given by the strategy the parents followed at time $n-1$ in environment $\omega_{n-1}$, that is $p_{\omega_{n-1}}$. This is the expressed in the following Lemma:
\begin{lem}\label{lem2}
For any $n\geq 1$, the population structure is conditionally independent of the past population process given the size of the population and the environment, i.e.
$$Z_{n}\Big |(|Z_{n}|,\omega_{n-1}) \quad \coprod \quad (Z_{0},\cdots,Z_{n-1})\Big |(|Z_{n}|,\omega_{n-1}).$$ 
\end{lem}
\noindent 
{\it Proof.}\ It is easily seen from the assumptions on the model \eqref{model} and from the non-hereditary assumption $\pi=\pi(\bar p)$ that the distribution of $Z_{n}$  given $(Z_0,\cdots Z_{n-1})$, $|Z_{n}|$ and $\omega_{n-1}$ is equal to the distribution of $\sum_{i=1}^{|Z_{n}|} 1_{\tau_{i}}$ with $\tau_i$ an iid sequence with distribution $p_{\omega_{n-1}}$. This distribution does not depend on $(Z_0,\cdots,Z_{n-1})$, this proves the conditional independence. \CQFD
\\ \ \\

The Lemma implies that given the environmental sequence $\omega$, the distribution of the population process $(Z_n)_{n\geq 1}$ is easily recovered from the size process $(|Z_n|)_{n\geq 1}$. Once again, this latter process turns out to be a simple BPRE and this allows us to compute the performance $\gamma(\bar p)$ of the strategy $\bar p$. Let $\omega^{(2)}=((\omega_{n-1},\omega_n))_{n\geq 1}$ denote the pair-environment with values in $\cE^2$. It is also stationary and ergodic and we denote by $\nu_2$ its stationary distribution, which  is the distribution of the pair $(\omega_1,\omega_2)$.

\begin{pr}\label{prop:pr5}
The size process $(|Z_n|)_{n\geq 1}$ is a simple branching process in environment $\omega^{(2)}$ with offspring distribution
$$ \Upsilon_{\bar p,(e_1,e_2)}= \int_{\cT} \Upsilon_{t,e_2}p_{e_1}(dt).$$
Conditionaly to $\omega$, the expected population size at time $n$ is 
$$\bbE_\omega[|Z_n|]=\bbE_{\omega_0}[|Z_1|]\prod_{k=1}^{n-1} m_{p_{\omega_{k-1},\omega_k}}$$
with $ m_{p_{e_1},e_2}=\int_{\cT} m_{t,e_2}p_{e_1}(dt) $ the first moment of $\Upsilon_{\bar p,(e_1,e_2)}$. Suppose the following integral exists, 
$$\gamma(\bar p)=\int_{\cE^2}\log\left( m_{p_{e_1},e_2}\right)\nu_2(de_1,de_2),$$ 
then 
$$\lim_{n\to\infty}n^{-1}\log \bbE_\omega[|Z_n|] =\gamma(\bar p)\quad \mbox{ almost\ surely.}$$
\end{pr}
\noindent
{\it Proof.}\ According to Lemma \ref{lem2}, given $(|Z_1|,\cdots,|Z_n|)$ and $\omega_{n-1}$, the population $Z_n$ has the same distribution as   
$\sum_{i=1}^{|Z_{n}|} 1_{\tau_{i}}$ with $\tau_i$ an iid sequence with distribution $p_{\omega_{n-1}}$. Intuitively, the $i$-th individual  has type $\tau_i$ chosen randomly on $\cT$ with distribution $p_{\omega_{n-1}}$. The size of the next generation $|Z_{n+1}|$ is then $\sum_{i=1}^{|Z_n|}\xi_{n,i}^{\tau_i,\omega_n}$ where $\xi_{n,i}^{\tau_i,\omega_n}$ is the offspring of the $i$-th individual of type $\tau_i$ in environment $\omega_n$. From this two step procedure i) random choice of the traits $t$ according to $p_{\omega_{n-1}}$, ii) reproduction with random offspring in environment $\omega_n$, we obtain the effective offspring distribution $\Upsilon_{p_{\omega_{n-1}},\omega_n}$ in environment $ \omega^{(2)}_n=(\omega_{n-1},\omega_n)$: it is the mixture of the offspring distributions $\Upsilon_{t,\omega_{n}}$, with mixing distribution $p_{\omega_{n-1}}$. Other properties are proved as in Proposition \ref{prop:pr1}.\CQFD
\\ \ \\
Note that the classification Theorem \ref{theo1} applies in this case as well and gives criteria for extinction or explosion of the population.

\subsection{Optimal strategies}
We now describe the set of optimal strategies when sensing is allowed. Is is also interesting to evaluate the gain between optimal strategy with or without sensing mechanisms. Optimality when sensing mechanism are allowed will be denoted with a double asterix whereas we keep a simple asterix for optimality without sensing. Let $\gamma^{\ast\ast}=\sup\{\gamma(\bar p);\ \bar p\in \cP(\cT)^{\cE}\}$ be the optimal growth rate when sensing is allowed and $\cP^{\ast\ast}=\{\bar p\in\cP(\cT)^{\cE};\ \gamma(\bar p)=\gamma^{\ast\ast}\}$ be the set of optimal strategies. 
Let $\nu_{e_1}(de_2)$ be the conditional distribution of $\omega_2$ given $\omega_1=e_1$, so that $\nu_2(de_1,de_2)=\nu_1(de_1)\nu_{e_1}(de_2)$. Note that conditional distributions are well-defined since $\cE$ is assumed to be a Polish space. From the previous section, the set of optimal strategies without sensing is denoted by $\cP^\ast$. It depends implicitly on the environment distribution $\nu_1$ and we write $\cP^\ast(\nu_1)$ to emphasize this dependence.  Suppose assumptions of Proposition \ref{prop:pr2} hold. Then optimal strategies with sensing are related with optimal strategies without sensing in the following way:

\begin{pr}\label{prop:pr6}
A strategy $\bar p$ is optimal if and only  if  
$$p_{e_1}\in \cP^\ast(\nu_{e_1})\quad \nu_1(de_1) {\rm\ almost\ everywhere.}$$
If conditions $({\bf C}_1)-({\bf C}_2)$ hold, then optimal strategies exists and form a closed convex set i.e.  
$\cP^{\ast\ast}$ is nonempty closed and convex.\\
If the family $\cM=\{e\mapsto m_{t,e}; t\in \cT\}$ is linearly independent, there is a unique optimal strategy $\bar p^{\ast\ast}$.
\end{pr}
\noindent
{\it Proof.} Using the explicit formula for $\gamma(\bar p)$ given in Proposition \ref{prop:pr5} and conditional probabilities, we compute
\begin{eqnarray*}
 \gamma(\bar p)&=&\int_{\cE^2}\log\left( m_{p_{e_1},e_2}\right)\nu_2(de_1,de_2)\\
&=& \int_{\cE}\nu_{1}(de_1)\int_{\cE} \log\left( m_{p_{e_1},e_2}\right)\nu_{e_1}(de_2)\\
&=& \int_{\cE}\gamma(p_{e_1},\nu_{e_1})\nu_{1}(de_1)
\end{eqnarray*}
whith $\gamma(p_{e_1},\nu_{e_1})=\int_{\cE} \log\left( m_{p_{e_1},e_2}\right)\nu_{e_1}(de_2)$ the growth rate associated with strategy $p_{e_1}$ in environment $\nu_{e_1}$. Hence $\gamma(\bar p)$ is optimal if we choose $p_{e_1}$ so that $\gamma(p_{e_1},\nu_{e_1})$ is maximal, i.e.  $p_{e_1}\in \cP^\ast(\nu_{e_1})$ almost surely. The properties of $\cP^{\ast\ast}$ and the uniqueness follow from the similar results for $\cP^\ast$ (cf Proposition \ref{prop:pr2}).\CQFD \\ \ \\

Then, the characterization of optimal strategies in the no-sensing case given in Proposition \ref{prop:pr3} directly extends to strategies with sensing as follows:
\begin{pr}\label{prop:pr6b} A strategy with sensing $\bar p\in\cP(\cT)^{\cE}$ is optimal if and only if 
$$ \int_{\cE} \frac{m_{t,e_2}}{m_{p_{e_1},e_2}}\nu_{e_1}(de_2)\leq 1,\quad {\rm for\ all\ } t\in\cT {\rm\ and\ }\nu_1(de_1)-{\rm almost\ everywhere}.$$
\end{pr}
{\it Proof}. The characterization follows directly from Proposition \ref{prop:pr3} and Proposition \ref{prop:pr6} together.\\

An interesting corollary states that no gain has to be expected from sensing mechanisms if the environment has some independence property. More precisely, 
\begin{cor}
Suppose that $\omega_1$ and $\omega_2$ are independent (i.e. $\nu_{2}=\nu_1\otimes \nu_1$), then the optimal growth rate with or without sensing are equal, i.e. $\gamma^\ast=\gamma^{\ast\ast}$.
\end{cor}
\noindent 
{\it Proof.} In this product case, the conditional distribution are trivial, i.e. $\nu_{e_1}\equiv \nu_1$ almost surely. Hence an optimal strategy is such that $p^{\ast\ast}_{e_1}\in\cP^\ast(\nu_1)$ and hence $\gamma(p^{\ast\ast}_{e_1},\nu_1)=\gamma^\ast$. Integrating with respect to $\nu_1(de_1)$, we obtain
$$\gamma^{\ast\ast}=\int_{\cE} \gamma(p^{\ast\ast}_{e_1},\nu_1) \nu_1(de_1)=\gamma^\ast.$$
This proves the result.
\CQFD \\ \ \\

Note that this result is rather intuitive: the sensing mechanism gives to the individual some information about the current environment state, but from the independence property, this is not useful for inference to the future environment state; hence the information is useless to decide which traits will be well-fitted in the next environment.

\subsection{Example: Finite dimensional case continued}
The results for optimal strategies without sensing developed in Section \ref{sec:fidi} together with Proposition \ref{prop:pr6} allow us to easily deduce the following properties for optimal strategies with sensing in the finite dimensional case.
\begin{pr}\label{prop:fidi2}
Suppose that ${\rm span}\cM$ has dimension $d$ and that conditions $({\bf C_1})-({\bf C_3})$ and $({\bf H})$ hold. Then 
there exists a unique optimal strategy with sensing $\bar p^{\ast\ast}\in\cP(\cT)^{\cE}$ such that
$$p^{\ast\ast}(e)=p^\ast(\nu_e) $$
where $p^\ast(\nu_e)$ is the optimal strategy without sensing from associated when the environment has marginal distribution $\nu_e$.
\end{pr}
Recall furthermore from Proposition \ref{prop:fidi1} that $p^\ast(\nu_e)$ is a discrete probability measure on $\cT$ with at most $d$ extremal types (that may depend on $e$). See Section \ref{s.non-hereditary/sensing} for more examples.

\subsection{Exemple: Haccou and Iwasa's example continued}
This is the sequel of subsection 3.5, the example by Haccau and Iwasa. Recall that the environment is given by a Gaussian stationary ergodic sequence $\omega=(\omega_n)_{n\geq 0}$ with stationary distribution $\nu_1=\cN(\mu,\sigma_2^2)$. Let $\rho\in (-1,1)$ be the pair correlation $\rho={\rm corr}(\omega_0,\omega_1)$. The case $\rho=0$ corresponds to independent environments, i.e. $\nu_2=\nu_1\otimes\nu_1$. Otherwise dependence holds and standard Gaussian computations give the conditional distribution $\nu_{e_1}=\cN(\mu+\rho(e_1-\mu),(1-\rho^2)\sigma_2^2)$. A particular realization of such a sequence is the Ornstein-Uhlenbeck sequence defined by 
$$\left\{\begin{array}{lll} 
\omega_0&=&\mu+\sigma_2 \eta_0 \\ 
\omega_{n}&=&\mu+\rho(\omega_{n-1}-\mu)+\sqrt{1-\rho^2}\sigma_2 \eta_{n},\quad n\geq 1 \\ 
\end{array}\right., $$
for i.i.d. standard normal innovations $(\eta_n)_{n\geq 0}$ (Gaussian white noise).

We have seen in section 3.5 that the optimal strategy without sensing is Gaussian (possibly degenerated) and we have given explicit formulas  for the parameters  \eqref{eq:optstrat} and for the corresponding growth rate \eqref{eq:optgrow}. According to Proposition \ref{prop:pr6}, we  deduce the optimal strategy when sensing is allowed: 
\begin{pr}\label{prop:HaccouIsa2} There exists a unique optimal strategy with sensing, which is denoted by $\bar p^{\ast\ast}$ and satisfies  $\nu_1(de)$-a.e. 
\begin{equation}\label{eq:optstrat2}
p^{\ast\ast}_{e}=\left\{\begin{array}{ll} \delta_{\mu+\rho(e-\mu)} & {\rm\ if\ } (1-\rho^2)\sigma_2^2\leq \sigma_1^2 \\ \cN(\mu+\rho(e-\mu),(1-\rho^2)\sigma_2^2-\sigma_1^2) & {\rm\ if\ }  (1-\rho^2)\sigma_2^2\geq \sigma_1^2 \end{array}\right..
\end{equation}
The corresponding optimal growth rate is given by
$$\gamma^{\ast\ast}= \left\{\begin{array}{ll} 
\log C-\frac{1}{2}\log(2\pi\sigma_1^2)-\frac{(1-\rho^2)\sigma_2^2}{2\sigma_1^2} & {\rm\ if\ } (1-\rho^2)\sigma_2^2\leq \sigma_1^2 \\ 
\log C-\frac{1}{2}\log(2\pi(1-\rho^2)\sigma_2^2)-\frac{1}{2} & {\rm\ if\ }  (1-\rho^2)\sigma_2^2\geq \sigma_1^2 
\end{array}\right. $$
\end{pr}
Finally, we can evaluate the relative gain of strategies with sensing over strategies without sensing: it gives an indication of the benefit that can be expected from sensing mechanisms. Let $\chi=\sigma_2^2/\sigma_1^2$. Then
$$\gamma^{\ast\ast}-\gamma^\ast=\left\{\begin{array}{ll} 
\frac{1}{2}\rho^2 \chi & {\rm\ if\ }  \chi\leq 1\\
\frac{1}{2}\log\chi-\frac{1}{2}(1-\rho^2)\chi +\frac{1}{2}& {\rm\ if\ } 1\leq \chi\leq (1-\rho^2)^{-1} \\ 
-\frac{1}{2}\log(1-\rho^2) & {\rm\ if\ }  \chi\geq (1-\rho^2)^{-1}
\end{array}\right.. $$
It is worth noting that this is an increasing function of the square correlation $\rho^2$: this indicates that the more correlated the random environment is, the more useful sensing mechanisms are. The intuitive idea is that higher correlations allows for more accurate prevision for the next environment and hence for a better fitted  offspring trait distribution in the environment to come.

\section{The hereditary case}\label{sec:nonher}
In the hereditary case, the trait distribution $\pi_{t,e}\in\cP(\cT)$ might depends on the trait of the parents. This dependency makes the study of the Lyapounov exponent $\gamma=\gamma(\pi)$ much more difficult because no reduction to a simple branching process in random environment is available. We have no explicit formula for $\gamma(\pi)$ in this case and determining the optimal strategy $\pi^\ast$ and the optimal growth rate $\gamma(\pi^\ast)$ might be very challenging. 

Nevertheless, we propose an interesting representation of the "finite time" growth rate
$$\gamma_n(\omega,\pi_0)=n^{-1}\log \bbE_{\omega,\pi_0}[|Z_n|]$$
in environment $\omega$ and  initial population consisting of a single individual with random trait with distribution $\pi_0$. The representation is in term of a functional of the $\cT$-valued Markov chain $T=(T_n)_{n\geq 0}$ in environment $\omega$ such that:\\
\begin{equation}\label{eq:defT}
\left\{\begin{array}{l}
\bbP_{\pi_0,\omega}(T_0\in \,\cdot\,)=\pi_0(\,\cdot\,)\\
\bbP_{\pi_0,\omega}(T_k\in \,\cdot\,|T_0=t_0,\cdots,T_{k-1}=t_{k-1})=\pi_{t_{k-1},\omega_{k-1}}(\,\cdot\,)\quad,\quad 1\leq k\leq n\\
\end{array}
\right..
\end{equation}
The Markov chain $T_n$ in random environment $\omega$ is time-heterogeneous because the transitions depend on time $n$ through the value of the environment $\omega_n$. However, in the no-sensing case when the transitions $\pi_{t,e}\equiv \pi_t$ do not depend on $e$, the Markov chain $T$ is time-homogeneous. The result is the following:
\begin{pr}\label{pr7} The finite time growth rate in environment $\omega$ and initial population consisting in a single individual with trait distributed according to $\pi_0$ is given by
$$\gamma_n(\omega,\pi_0)= n^{-1}\log \bbE_{\omega,\pi_0}[M_n(T,\omega)]$$
with 
$$M_n(T,\omega)=\prod_{k=0}^{n-1}m_{T_k,\omega_k}.$$
\end{pr}
For the sake of clarity and conciseness, this proposition will be proved together with Theorem \ref{theo2} below.  
To our best knowledge, there is no simple way  to deal with the asymptotics behaviour of $\gamma_n$ and this issue might be challenging. We provide a lower bound for the growth rate $\gamma_n(\omega,\pi_0)$ that might be more tractable. 
The lower bound for $\gamma_n(\omega,\pi_0)$ is obtained using Jensen's inequality (with the concave function $\log$):
$$\gamma_n(\omega,\pi_0)= n^{-1}\log \bbE_{\omega,\pi_0}[M_n(T,\omega)]\geq n^{-1}\bbE_{\omega,\pi_0}[\log M_n(T,\omega)].$$
The advantage here is that the lower bound 
$$n^{-1}\bbE_{\omega,\pi_0}[ \log M_n(T,\omega)]=n^{-1}\sum_{k=0}^{n-1}\bbE_{\omega,\pi_0}[\log m_{T_k,\omega_k}] $$
can be written as an additive functional of the environment $\omega$ and we can then use the ergodic theorem to control the convergence. Note that this technique can be refined using changes of measures. Let $\tilde \pi_{t,e}$ any kernel family and denote by $\tilde T$ the $\cT$-valued Markov chain in environment $\omega$ starting from distribution $\pi_0$ and with transitions given by $\tilde\pi_{t,e}$ (given by equations similar to \eqref{eq:defT}). We suppose that $\pi_{t,e}$ is absolutely continuous with respect to $\tilde\pi_{t,e}$,  i.e. $\pi_{t,e}(dt')=f_{t,e}(t')\tilde\pi_{t,e}(dt')$ where $f_{t,e}$ stands for the density of $\pi_{t,e}$ with respect to $\tilde\pi_{t,e}$. Then, using changes of measures, we have
$$\bbE_{\omega,\pi_0}[M_n(T,\omega)]=\bbE_{\omega,\pi_0}[\tilde M_n(\tilde T,\omega)] $$ 
with 
$$ \tilde M_n(\tilde T,e)=\prod_{k=0}^{n-1}m_{\tilde T_k,\omega_k}f_{\tilde T_k,\omega_k}(\tilde T_{k+1}).$$
Using this, the lower bound becomes
$$\gamma_n(\omega,\pi_0)=n^{-1}\bbE_{\omega,\pi_0}[\log \tilde M_n(\tilde T,\omega)]\geq n^{-1}\sum_{k=0}^{n-1}\bbE_{\omega,\pi_0}\left[\log\left( m_{\tilde T_k,\omega_k}f_{\tilde T_k,\omega_k}(\tilde T_{k+1})\right)\right]. $$

\section{Typical genealogies: a mean field approach}
\subsection{Convergence of the typical genealogy in the infinite population limit}
As explained in Baake and Georgii \cite{BaaGeo1,BaaGeo2}, the evolution of a branching population can be studied from two possible perspectives: either forward or backward in time. So far, we have focused on the first point of view and mainly studied the growth rate of the population after a large numbers of generations. By way of contrast, the backwards or retrospective aspect of the population concerns the lineages extending back into past from the presently living individuals and asks for the characteristics of the ancestors along such lineages. We now turn to this second perspective and wonder what is the typical lineage or genealogy (backward in time) of an individual chosen at random in the $n$-th generation.

Some definitions are needed here. The right formalism to keep track of the genealogy is the formalism of labeled rooted trees and forests, where the trees stands for the descendence of each ancestor represented by a root, and labels keep track of the traits of the individuals. However we keep this formalism to its minimum. Let $\cG_n$ denote the population at the $n$-th generation. To each individual $g\in \cG_n$, we associate its lineage or genealogy $\ell(g)=(t_0,\cdots,t_n)\in\cT^{n+1}$ with the interpretation that $t_n$ is the trait of $g$, and $t_{k}$ the trait of his ancestor in the $k$-th generation $\cG_k$, $0\leq k\leq n-1$. The typical genealogy is defined as the genealogy of an individual chosen at random in the $n$-th generation. This obviously requires the $n$-generation to be non empty, in which case we adopt the convention that the typical genealogy is $\emptyset$. The distribution of the typical genealogy is given by
$$\pi_{n}=\left\{\begin{array}{ll}
\frac{1}{{\rm card} \cG_n}\sum_{g\in\cG_n}\delta_{\ell(g)}&{\rm\ if \ }\cG_n\neq\emptyset \\
\delta_{\emptyset}&{\rm\ if \ }\cG_n=\emptyset \\
\end{array}\right..$$
This is a random measure on $\cT^{n+1}\cup\{\emptyset\}$. We denote by $\bbP_{N,\pi_0,\omega}$ the probability measure corresponding to a population evolving in environment $\omega$, starting a time $0$ from $N$ individuals with traits i.i.d. with distribution $\pi_0$. In the following, we focus on the typical genealogy in the infinite population limit $N\to\infty$. Recall the definition of the Markov chain $T=(T_n)_{n\geq 0}$ given by equation \eqref{eq:defT}. We have the following mean field result:

\begin{theo}\label{theo2}
Under the probability $\bbP_{N,\pi_0,\omega}$, the typical genealogy distribution $\pi_{n}$ almost surely weakly converges as $N\to\infty$ to the distribution $\hat \pi_{n,\pi_0,\omega}$ defined by 
$$\hat \pi_{n,\pi_0,\omega}(A)=\frac{\bbE_{\omega,\pi_0}[M_n(T,\omega)1_{(T_0,\cdots,T_n)\in A}]}{\bbE_{\omega,\pi_0}[M_n(T,\omega)]}\quad,\quad A\subset \cT^{n+1}.$$
\end{theo}
\ \\ \ \\
Recall from Proposition \ref{pr7} that ${\bbE_{\omega,\pi_0}[M_n(T,\omega)]}=\exp(n\gamma_n(\pi_0,\omega))$ is the mean number of individuals in the $n$-th generation of a population evolving in environment $\omega$ and starting from a single individual with random trait with distribution $\pi_0$. 
\\ \ \\
\noindent
{\it Proof of Theorem \ref{theo2} and Proposition \ref{pr7}.} Let $A=A_0\times\cdots\times A_n$ be a product subset of $\cT^{n+1}$. The number of individuals in the $n$-th generation with genealogy in $A$ is 
$$\cN_n(A)=\sum_{g\in\cG_n}\delta_{\ell(g)}(A)$$
and $\cN_n(\cT^{n+1})$ denotes the total number of individuals in $\cG_n$. 
We can see $\cN_n$ as a non-normalized measure and $\pi_n$ is the probability measure associated with by the relation
$$\pi_{n}=\frac{1}{\cN_n(\cT^{n+1})}1_{\{\cN_n(\cT^{n+1})>0\}}\cN_n+1_{\{\cN_n(\cT^{n+1})=0\}}\delta_{\emptyset}.$$
From the branching property, the distribution of $\cN_n(A)$ under $\bbP_{N,\pi_0,\omega}$ is equal to the sum of $N$ independent copies $\sum_{i=1}^N \cN_n^{(i)}(A)$ under  $\bbP_{1,\pi_0,\omega}^{\otimes N}$. As a consequence of the weak law of large numbers, the distribution of $\frac{1}{N}\cN_n(A)$ under $\bbP_{N,\pi_0,\omega}$ weakly converge to $\bbE_{1,\pi_0,\omega}(\cN_n(A))$. The results also holds for $A=\cT^{n+1}$ and taking the quotient, we see that
under $\bbP_{N,\pi_0,\omega}$, $\pi_{n}(A)$ weakly converge to 
$$\hat \pi_{n,\pi_0,\omega}(A)=\frac{\bbE_{1,\pi_0,\omega}(\cN_n(A))}{\bbE_{1,\pi_0,\omega}(\cN_n(\cT^{n+1}))}$$
provided the denominator is non zero. Note that 
$$\gamma_n(\omega,\pi_0)=n^{-1}\log \bbE_{1,\pi_0,\omega}(\cN_n(\cT^{n+1})).$$
The theorem and the proposition (letting $A=\cT^{n}$) are then a consequence of 
\begin{eqnarray*}
\bbE_{1,\pi_0,\omega}(\cN_n(A))&=&\int_{A}\pi_0(dt_0) \prod_{k=0}^{n-1} m_{t_k,\omega_k}\pi_{t_k,\omega_k}(dt_{k+1})\\
&=& \bbE_{\omega,\pi_0}[M_n(T,\omega)1_{(T_0,\cdots,T_n)\in A}].
\end{eqnarray*}
This last relation is proven by induction: individuals in the $n$-th generation with genealogy in $A$ are the offsprings with traits in $A_n$ of individuals in the $n-1$ generation with genealogy in $A_0\times\cdots\times A_{n-1}$. This yields
$$ \bbE_{1,\pi_0,\omega}(\cN_n(A))= \int_{A_{n-1}\times A_n}\bbE_{1,\pi_0,\omega}(\cN_{n-1}(A_0\times\cdots\times A_{n-2}\times dt_{n-1}))m_{t_{n-1},\omega_{n-1}}\pi_{t_{n-1},\omega_{n-1}}(dt_n).$$
For $n=0$, $\bbE_{1,\pi_0,\omega}(\cN_0(A_0))=\pi_0(A)$. It remains to note that $\bbE_{1,\pi_0,\omega}(\cN_n(\cT^{n+1})) \neq 0$ and this implies that the mass of $\emptyset$ vanishes in the limit.\CQFD
\\ \ \\

\subsection{The typical genealogy in the non-hereditary case}
When the traits are non-hereditary, the trait distributions  $\pi_{t,e}$ do not depend on $t$ and the mean field typical genealogy distribution $\hat\pi_{n,\omega,\pi_0}$ has a very simple form. The following Proposition is given in the context of a population with a sensing mechanism $\pi=\pi((p_e)_{e\in\cE})$. The no-sensing case corresponds to the particular case when $p_e\equiv p$. 
\begin{pr}
In the non-hereditary case, the typical genealogy distribution $\hat\pi_{n,\omega,\pi_0}$ is the product measure on $\cT^{n+1}$ defined by
$$\hat\pi_{n,\omega,\pi_0}(dt_0,\cdots,dt_n)=\otimes_{i=0}^n \hat \pi_{i,\omega_i}(dt_i) $$
where
$$\left\{\begin{array}{lll} 
\hat\pi_{0,\omega_0}(dt_0)&=&\frac{m_{t_0,\omega_0}}{m_{\pi_0,\omega_0}}\pi_0(dt_0),\\
\hat \pi_{i,\omega_i}(dt_i) &=&\frac{m_{t_i,\omega_i}}{m_{p_{\omega_{i-1}},\omega_i}}p_{\omega_{i-1}}(dt_i), \quad 1\leq i\leq n-1,\\
\hat \pi_{n,\omega_n}(dt_n)&=&p_{\omega_{n-1}}(dt_n)\end{array}\right..$$
\end{pr}
{\it Proof.} In the non-hereditary case, the Markov chain in random environment $\omega$ defined by equation \eqref{eq:defT} is simple because conditionally on the environment $\omega$, the random variables $T_0,\cdots,T_n$ are independent with $T_0$  distributed as  $\pi_0$ and, for $1\leq i\leq n$, $T_i$ distributed according to $p_{\omega_{i-1}}$. Let $A=A_0\times\cdots\times A_n\subset \cT^{n+1}$. Using independence and Theorem \ref{theo2}, we compute: \begin{eqnarray*}
\hat \pi_{n,\pi_0,\omega}(A)&=&\frac{\bbE_{\omega,\pi_0}[M_n(T,\omega)1_{(T_0,\cdots,T_n)\in A}]}{\bbE_{\omega,\pi_0}[M_n(T,\omega)]}\\
&=&\frac{\bbE_{\omega,\pi_0}[m_{T_0,\omega_0}1_{T_0\in A_0}]\times\cdots\times \bbE_{\omega,\pi_0}[m_{T_{n-1},\omega_{n-1}}1_{T_{n-1}\in A_{n-1}}]\bbE_{\omega,\pi_0}[1_{T_{n}\in A_{n}}] } {\bbE_{\omega,\pi_0}[m_{T_0,\omega_0}]\times\cdots\times \bbE_{\omega,\pi_0}[m_{T_{n-1},\omega_{n-1}}]}.
\end{eqnarray*} 
This proves the independence property and gives the marginal distributions: 
$$\hat \pi_i(A_i)=\frac{\bbE_{\omega,\pi_0}[m_{T_{i},\omega_{i}}1_{T_{i}\in A_{i}}]}{\bbE_{\omega,\pi_0}[m_{T_{i},\omega_{i}}]},\quad 1\leq i\leq n-1.$$
\CQFD 
\\ \ \\
The interpretation of the above proposition is the following: the mean field typical genealogy in environment $\omega$ consists in independent traits $\hat T_0,\cdots, T_n$ where\\
- the distribution of $\hat T_0$ is a biased version of $\pi_0(dt_0)$ with bias function equal to $m_{t_0,\omega_0}$ the mean number of offspring's of an individual of type $t_0$,\\
- the distribution of $\hat T_i$ is a biased version of $p_{\omega_i}(dt_i)$ with bias function equal to $m_{t_i,\omega_i}$,\\
- the distribution of $\hat T_n$ is  $p_{\omega_i}(dt_i)$ (there is no bias because the offspring of the last generation $\cG_n$ is not involved since the population is considered until time $n$ only).

\subsection{The typical genealogy in a hereditary case: Haccou and Isawa's example continued}
Recall that $\cT=\cE=\bbR$, $\omega$ is a Gaussian stationary ergodic sequence with stationary distribution $\cN(\mu,\sigma_2^2)$ and $m_{t,e}$ is given by 
$$m_{t,e}=\frac{C}{\sqrt{2\pi\sigma_1^2}}\exp\left(-\frac{(t-e)^2}{2\sigma_1^2}\right).$$
It remains to precise the trait distributions $\pi_{t,e}$ and $\pi_0$. To make the model explicitly solvable, we require the Markov chain $ T=(T_0,\cdots,T_n)$ in environment $\omega$  to be multivariate Gaussian. This impose that the transition $\pi_{t,e}$ are of the form $\pi_{t,e}=\cN(\alpha_e t+\beta_e,\theta_e)$ with 
$$\bbE[T_1|\omega_0=e,T_0=t]=\alpha_e+\beta_e t \quad {\rm and}\quad {\rm Var}[T_1|\omega_0=e,T_0=t]=\theta_e^2.$$
Let $\pi_0=\cN(\mu_0,s^2_0)$ be the initial distribution. Alternatively, we have the representation in environment $\omega$
$$\left\{ \begin{array}{lll}
T_0&=&\mu_0+s_0N_0, \\
T_{k+1}&=& \alpha_{\omega_k} +\beta_{\omega_k}T_k+\theta_{\omega_k}N_{k+1}\quad 0\leq k\leq n-1,
\end{array}\right.
$$
with $N_0,\cdots,N_n$ independent standard normal variable. We introduce the $(n+1)\times 1$ vectors
$$T=\left(\begin{array}{c}T_0\\ T_1\\ \vdots \\ T_{n} \end{array}\right)\quad,\quad  N=\left(\begin{array}{c}N_0\\ N_1\\ \vdots \\ N_{n} \end{array}\right)\quad,\quad A_\omega=\left(\begin{array}{c}\mu_0\\ \alpha_{\omega_0}\\ \vdots \\ \alpha_{\omega_{n-1}} \end{array}\right)$$
and the $(n+1)\times (n+1)$ matrices 
$$C_\omega=\left(\begin{array}{ccccc}0 & 0 & \cdots  & 0 & 0 \\ \beta_{\omega_0} & 0 & \cdots & 0 & 0 \\ 0 & \beta_{\omega_1} & 0 &  0 & 0\\ 0 & 0 & \ddots & 0 & 0 \\0 & 0  & \cdots & \beta_{\omega_{n-1}} & 0 \end{array}\right)\quad ,\quad S_\omega=\left(\begin{array}{ccccc}s_0^2 & 0 & \cdots  & 0 & 0 \\ 0 & \theta_{\omega_0}^2 & \cdots & 0 & 0 \\ 0 & 0 & \theta_{\omega_1}^2 &  0 & 0\\ 0 & 0 & \cdots & \ddots & 0 \\0 & 0  & \cdots & 0 & \theta_{\omega_{n-1}}^2 \end{array}\right). $$
With these notations, the recursive relation turns into
$$T=A_\omega+C_\omega T+S_\omega N $$
and this yields
$$T=(Id-C_\omega)^{-1}(A_\omega +S_\omega N). $$
We deduce the mean $\mu$ and covariance matrix $\Sigma$ for the Gaussian vector $T$
$$\mu=(Id-C_\omega)^{-1}A_\omega \quad {\rm and}\quad \Sigma= (Id-C_\omega)^{-1}S_\omega(Id-C_\omega^{'})^{-1}.$$
Suppose that $S_\omega$ is invertible, then $ T$ has density
$$f_{ T}( t)=(2\pi)^{-(n+1)/2}\det(\Sigma)^{-1/2}\exp\left(-\frac{1}{2}( t-\mu)'\Sigma^{-1}( t-\mu) \right).$$
Denote by $\hat T=(\hat T_0,\cdots,\hat T_n)$ a random vector with distribution $\hat \pi_{n,\omega,\pi_0}$. 
According to Theorem \ref{theo2}, $\hat T$ is a biased version of $T$ and has a density given by 
\begin{eqnarray*}
f_{\hat T}( t)& = &  e^{-n\gamma_n(\omega,\pi_0)}M_n( t,\omega)f_{T}( t)\\
&=&e^{-n\gamma_n(\omega,\pi_0)} C^n\sigma_1^{-n}(2\pi)^{-n-1/2}\det(\Sigma)^{-1/2}\exp\left(-\frac{1}{2}(t-\mu)'\Sigma^{-1}(t-\mu)-\frac{1}{2\sigma_1^2}\sum_{i=0}^{n-1}(\omega_i-t_i)^2 \right).
\end{eqnarray*}
Introducing $J_{n,1}$ the diagonal matrix $J_{n,1}=(\delta_{0\leq i=j\leq n-1})_{0\leq i,j\leq n}$ and $V_{n,\omega}$ the column vector $V_{n,\omega}=(\omega_0,\cdots,\omega_{n-1},0)'$, the last exponential factor rewrites
$$\exp\left(-\frac{1}{2}t'(\Sigma^{-1}+\sigma_1^{-2}J_{n,1})t+(\Sigma^{-1}\mu+\sigma_1^{-2}V_{n,\omega})'t -\frac{1}{2}\mu'\Sigma^{-1}\mu-\frac{1}{2}\sigma_1^{-2}V_{n,\omega}'V_{n,\omega} \right).$$
We recognize that $f_{\hat T}$ is a multivariate Gaussian density with mean $\hat \mu$ and covariance matrix $\hat\Sigma$ of the form 
$$f_{\hat T}( t)=(2\pi)^{-(n+1)/2}\det(\hat \Sigma)^{-1/2}\exp\left(-\frac{1}{2}( t-\hat\mu)'\hat\Sigma^{-1}( t-\hat\mu) \right).$$
Identifying both expressions, we obtain after simplification
\begin{eqnarray}
\hat\Sigma&=&\left(\Sigma^{-1}+\sigma_1^{-2} J_{n,1}\right)^{-1}, \label{eq:defhatsigma}\\
\hat\mu&=&\hat\Sigma\left(\Sigma^{-1}\mu+\sigma_1^{-2}V_{n,\omega}\right)=\left(\Sigma^{-1}+\sigma_1^{-2} J_{n,1}\right)^{-1}\left(\Sigma^{-1}\mu+\sigma_1^{-2}V_{n,\omega}\right), \label{eq:defhatmu}\\
\gamma_n(\omega,\pi_0)&=&\log \frac{C}{\sqrt{2\pi}\sigma_1}-\frac{1}{2n}\log \det(I_{n+1}+\sigma_1^{-2}J_{n,1}\Sigma)+\frac{1}{2n}\left(\hat\mu'\hat\Sigma^{-1}\hat\mu-\mu'\Sigma\mu-\sigma_1^{-2}V_{n,\omega}'V_{n,\omega} \right).\label{eq:growrate}
\end{eqnarray}
These computations prove the following result:
\begin{pr}\label{prop8}
In the Haccou and Isawa model with Markov Gaussian environment, the typical genealogy distribution $\hat \pi_{n,\omega,\pi_0}$ is the Gaussian distribution on $\bbR^{n+1}$ with mean $\hat\mu $ given by  \eqref{eq:defhatmu} and covariance matrix $\hat\Sigma$ given by \eqref{eq:defhatmu}. Furthermore, the finite time growth rate $\gamma_n(\omega,\pi_0)$ is given by \eqref{eq:growrate}.
\end{pr}

\end{document}